\documentclass{article}
\usepackage{times}
\usepackage{fullpage}
\usepackage{latexsym,amsfonts,amsmath}
\usepackage{graphics,epic,eepic,color}
\usepackage{ifpdf}
\usepackage{algorithm,algorithmicx,algpseudocode}

\newtheorem{theorem}{Theorem}

\newtheorem{lemma}{Lemma}

\newenvironment{proof}{\noindent{\bf Proof:}\ }{\(\qed\) \par\medskip}

\newcommand{\qed}{\quad\mbox{\rule{7pt}{7pt}}}

\newcommand{\cmt}[1]{\ifhmode\newline\fi{\sf *** \ \ #1 \\}}

\newcommand{\NN}{\mathbf{N}}

\newcommand{\RR}{\mathbf{R}}
\newcommand{\ones}{\mathbf{1}}
\newcommand{\zeros}{\mathbf{0}}
\newcommand{\E}{\mathbf{E}}
\newcommand{\LL}{L}
\newcommand{\DD}{D}

\newcommand{\calE}{\mathcal{E}}
\newcommand{\calF}{\mathcal{F}}

\newcommand{\calU}{\mathcal{U}}

\newlength{\saveparindent}
\def\eproof{\end{rm}\addtolength{\parskip}{-4pt}%
\setlength{\parindent}{\saveparindent}}
\newcommand{\bprooff}[1]{\begin{rm}\protect\vspace{5pt}%
\noindent\textbf{Proof of #1: }\addtolength{\parskip}{4pt}%
\setlength{\parindent}{0pt}}
\newenvironment{prooff}[1]{\par\bprooff{#1}}{\eproof\(\qed\)\par}

\begin{document}

\title{On a Speculated Relation Between Chv\'atal-Sankoff
  Constants of Several Sequences}

\author{
  M.~Kiwi\thanks{
 Departamento de Ingenier\'{i}a ~Matem\'{a}tica. Centro de~Modelamiento Matem\'{a}tico (UMI 2807, CNRS), U.~Chile.
  Web: \texttt{www.dim.uchile.cl/$\sim$mkiwi}.
  Gratefully acknowledges the support of
    CONICYT via FONDAP in Applied Mathematics and
    Anillo en Redes ACT08.}
  \and
J.~Soto\thanks{Department of Mathematics,
  MIT.
   \texttt{jsoto@math.mit.edu} .
  Gratefully acknowledges the support
    of CONICYT via Anillo en Redes ACT08.}
}

\date{}

\maketitle

\begin{abstract}
It is well known that, when normalized by $n$, the expected length
  of a longest common subsequence of $d$ sequences of length $n$ over
  an alphabet of size $\sigma$
  converges to a constant $\gamma_{\sigma,d}$.
We disprove a speculation by Steele regarding a possible relation
  between $\gamma_{2,d}$ and $\gamma_{2,2}$. In order to do that we also obtain some new lower bounds for $\gamma_{\sigma,d}$, when both $\sigma$ and $d$ are small integers.
\end{abstract}

\section{Introduction}

String matching is one of the most intensively analyzed problems
  in computer science.
Among string matching problems the longest common subsequence
  problem (LCS) stands out.
This problem consists of finding the longest
  subsequence common to all strings in a set of sequences
  (often just two).
The LCS problem is the basis of Unix's \texttt{diff}
  command, has applications in bioinformatics, and also
  arises naturally in remarkably distinct
  domains such as cryptographic snooping, the mathematical
  analysis of bird songs, and comparative genomics.
In addition, the LCS problem offers a concrete basis for the illustration
  and benchmarking of mathematical methods and tools such as subadditive
  methods and martingale inequalities; see for example Steele's
  monograph~\cite{steele86}.

Although the LCS problem has been studied under many different contexts
  there are several issues concerning it that are still unresolved.
The most prominent of the outstanding questions
  relating to the LCS problem concerns the length $\LL_{n,\sigma,d}$ of a LCS
  of $d$ sequences of $n$ characters chosen uniformly
  and independently over some alphabet of size~$\sigma$.
Subadditivity arguments yield that for fixed $d$ and
  $n$ going to infinity,
  the expected value of $\LL_{n,\sigma,d}$ normalized by $n$ converges
  to a constant $\gamma_{\sigma,d}$.
For $d,\sigma\geq 2$, the precise value of
  $\gamma_{\sigma,d}$ is unknown.
The constant $\gamma_{2,2}$ is referred to as the Chv\'atal-Sankoff
  constant.
The calculation of its exact value is an over 3 decades old open problem.
The determination of its value has received a fair amount of
  attention, starting with the work of Chv\'atal and Sankoff~\cite{cs75},
  encompassing among
  others~\cite{deken79,alexander94,dp95,dancik98,bgns99,lueker03},
  and is explicitly stated in several well known texts such
  as the ones by
  Waterman~\cite[\S~11.1.3]{waterman95}, Steele~\cite[p.~3]{steele97},
  Pevzner~\cite[p.~107]{pevzner00},
  and~Szpankowski~\cite[p.~109]{szpankowski00}.
To the best of our knowledge the current sharpest bounds on $\gamma_{2,2}$
  are due to Lueker~\cite{lueker03}
  who established that $0.788071\leq \gamma_{2,2}\leq 0.826280$.

The starting point for this investigation is
  the following comment by Steele~\cite{steele86}:
\begin{quote}
``It would be of interest to relate
  $c_3$ to $c^2$, and one is tempted to speculate
  that $c_3=c^{2}$ (and more generally
  that $c_k=c^{k-1}$).
Computational evidence does not yet rule this out.''
\end{quote}
Here, Steele uses $c$ to denote the limiting value of the longest
common subsequence of two random sequences of length $n$ normalized by
$n$ as $n$ goes to infinity, and in general, he uses~$c_k$ to denote
the analogous constant for $k$ sequences. However, it is unclear if in
this comment he uses~$c$ and $c_k$ to denote the constants
$\gamma_{2,2}$ and $\gamma_{k,2}$ (i.e. specifically for the case of
alphabet size 2) or if he is generically denoting the constants for
arbitrary alphabet size.
Dan\v{c}\'{\i}k~\cite{dancik98} cites the previous statement
  as a conjecture by Steele using the second interpretation, i.e.,
  as the claim that for all $d\geq 3$ and $\sigma\geq 2$,
\begin{equation}\label{eq:dancik-conjecture}
\gamma_{\sigma,d}=\gamma_{\sigma,2}^{d-1}\,.
\end{equation}
Dan\v{c}\'{\i}k~\cite[Theorem~2.1, Corollary~2.1]{dancik98} shows that for $d\geq 2$
\begin{equation*}
  1 \leq \liminf_{\sigma\to\infty} \sigma^{1-1/d}\gamma_{\sigma,d}
  \leq
  \limsup_{\sigma\to\infty} \sigma^{1-1/d}\gamma_{\sigma,d}
  \leq e\,.
\end{equation*}
Hence, if~(\ref{eq:dancik-conjecture}) was true, then
  for $\epsilon>0$ and $\sigma$ sufficiently large,
\[
1-\epsilon \leq \sigma^{1-1/d}\gamma_{\sigma,d}
  = \sigma^{1-1/d}\gamma_{\sigma,2}^{d-1}
  \leq \sigma^{1-1/d}\left(\frac{e(1+\epsilon)}{\sqrt{\sigma}}\right)^{d-1}\,.
\]
Dan\v{c}\'{\i}k's results
  disprove~(\ref{eq:dancik-conjecture}) by observing
  that for $d>2$ one may choose $\sigma$ large enough so
  as to make the rightmost term of the last displayed equation
  arbitrarily close to $0$.

If we use the first interpretation of Steele's speculation quoted above,
  i.e., considering only the case of binary alphabets as we believe it
  was intended, then~(\ref{eq:dancik-conjecture})
  is not invalidated by Dan\v{c}\'{\i}k's work.

In~\cite{steele86}, Steele does not justify his speculation.
The following non-rigorous argument
  gives some indication that one should expect
  that $\gamma_{2,3}$ is strictly bigger than $\gamma_{2,2}^{2}$.
Indeed, let $A_1$, $A_2$ and $A_3$ be three
  independently and uniformly chosen binary sequences of length $n$.
For $i\neq j$ and very large values of $n$ one knows that a longest common
  subsequence $\ell_{i,j}$ of sequences $A_i$ and~$A_j$
  would be of length approximately $\gamma_{2,2}n$.
One would expect (although we can not prove it) that~$\ell_{i,j}$
  would behave like a uniformly chosen binary string of length
  $\gamma_{2,2}n$.
Sequences $\ell_{1,2}$ and~$\ell_{2,3}$ are clearly correlated.
However, one might guess that the correlation is weak (again, we
  can certainly neither formalize nor prove such a statement).
The previously stated discussion suggests that
  a longest common subsequence $\ell_{1,2,3}$ of $\ell_{1,2}$ and~$\ell_{2,3}$
  should be of length approximately $\gamma_{2,2}^{2}n$.
Since $\ell_{1,2,3}$ is clearly a longest common subsequence of
  $A_{1}$, $A_{2}$ and $A_{3}$, one is led to
  conclude that
\begin{equation}
\gamma_{2,3}\geq \gamma_{2,2}^{2}\,.\label{eq:main}
\end{equation}
However, there are two good reasons why one suspects that this last
  inequality should be strict:
\begin{itemize}
\item Since $\ell_{2,3}$ has only a fraction of
  $A_{3}$'s length, one expects that a longest common
  subsequence of $\ell_{1,2}$ and $A_{3}$ is significantly
  larger than a longest common subsequence of
  $\ell_{1,2}$ and $\ell_{2,3}$.
\item The longest common subsequence of $A_{1}$,
  $A_{2}$ and $A_{3}$ might arise by taking a longest
  common subsequence on sub-optimal common subsequences
  $\ell'_{1,2}$ and $\ell'_{2,3}$ of $A_{1}$ and $A_{2}$,
  and $A_{2}$ and $A_{3}$, respectively.
\end{itemize}
This work's main contribution is to show that
  the inequality in~(\ref{eq:main}) is indeed strict.

\medskip
In Section~\ref{sec:simple} we give a simple argument that proves
  that when $\sigma$ is fixed and $d$ is large
  the identity $\gamma_{\sigma,d}=\gamma_{\sigma,2}^{d-1}$
  does not hold.
The underlying argument is essentially an application of the probabilistic
  method.
However, it might still be possible that the relation would hold for
  some specific values of $\sigma$ and $d$.
Of particular interest is the case of binary sequences,
  i.e.~$\sigma=2$.
In Section~\ref{sec:complex} we show that even this weaker identity
  does not hold, i.e.~that $\gamma_{2,3}\neq \gamma_{2,2}^{2}$.
To achieve this goal, we rely on Lueker's~\cite{lueker03}
  $U=0.826280$  upper bound on $\gamma_{2,2}$ and determine
  a lower bound on $\gamma_{2,3}$ which is strictly larger
  than $U^{2}\geq \gamma_{2,2}^{2}$.
The lower bound on $\gamma_{2,3}$ is obtained by an approach similar to
  the one used by Lueker~\cite{lueker03} to lower bound $\gamma_{2,2}$,
  although in our case we have to consider a non-binary alphabet.
Aside from the extra notation needed to handle the cases $\sigma, d>2$, our
  treatment is a straightforward generalization of the approach used by
  Lueker.
(In fact, in order to keep the exposition as clear as possible
  we do not even use the optimization tweaks implemented by Lueker in order
  to take advantage of the symmetries inherent to the problem and
  objects that arise in its analysis.)
We conclude with some final comments in Section~\ref{sec:conclusion}.

\section{Disproving $\gamma_{\sigma,d}=\gamma_{\sigma,2}^{d-1}$
  for large $d$}\label{sec:simple}
We start this section by introducing some notation. Given strings
$A_1,\ldots, A_d$ of length $n$, we denote by
$L(A_1,\ldots,A_d)$ the length of the longest common subsequence of
all $A_{i}$'s. Let $\calU_{n,\sigma}$ be the distribution of sequences
of length $n$ whose characters are chosen uniformly and independently
from $\Sigma = \{1,\ldots, \sigma\}$. We denote by
$\LL_{n,\sigma,d}$ the random variable $L(A_1,\ldots,A_d)$ when all
the $A_i$ are chosen according to $\calU_{n,\sigma}$. Finally, we let
$\gamma_{\sigma,d}$ denote the limit of $\E\LL_{n,\sigma,d}/n$ when
  $n\to\infty$
(the existence of this limit follows from standard subadditivity
arguments~\cite{cs75}).

In what follows, we give a lower bound for $\gamma_{\sigma,d}$ that is
  independent of $d$.
This bound is based on the following simple fact: If $X$ is chosen
  according to $\calU_{n,\sigma}$ and $n$ is large, then the number of
  occurrences of a fixed character in~$\Sigma$ is roughly
  $n/\sigma$.
Intuitively, this means that for a set of $d$ random strings of
  (very large) length $n$, with very high probability
  a sequence formed by roughly~$\lfloor n/\sigma\rfloor$
  equal characters will be
  a common subsequence of all the $d$ random strings.

\begin{lemma}\label{lem:simple}
For all $d$ and $\sigma$, we have
  $\gamma_{\sigma,d} \geq 1/\sigma.$
\end{lemma}
\begin{proof}
Let $A_1,\ldots, A_d$ be $d$ independent random strings chosen
  according to $\calU_{n,\sigma}$.
Let $X_i$ denote the number of times the character $c\in\Sigma$ appears in
  $A_i$, and $X=\min\{X_1,\ldots, X_d\}$.
The string $c^X$ formed by $X$ copies of the
  character $c$ is a common subsequence of all $X_{i}$'s.
It follows that $L(A_1,\ldots,A_d) \geq X$.

Each $X_i$ is a binomial variable with parameter $p=1/\sigma$.
By a standard Chernoff bound~\cite[Remark~2.5]{jlr00} we have that for any
  $0 < \varepsilon <1$,
\[
\Pr[X_i  \leq (1-\varepsilon)np ] \leq \exp(-2n(p\varepsilon)^2).
\]
Applying Markov's inequality, and recalling that the $X_i$'s are
  independent, it follows that:
\[
\E X \geq (1-\varepsilon)np\Pr[X \geq (1-\varepsilon)np]\\
      \geq (1-\varepsilon)np[1- \exp(-2n(p\varepsilon)^2)]^d.
\]
Letting $n$ be sufficiently large so that $[1- \exp(-2n(p\varepsilon)^2)]^d \geq (1-2\varepsilon)/(1-\varepsilon)$, we obtain $\E X \geq np(1-2\varepsilon)$. Therefore:
  \[\frac{\E \LL_{n,\sigma,d}}{n} = \frac{\E L(A_1,\ldots, A_d)}{n} \geq \frac{\E X}{n} \geq (1-2\varepsilon)p = \frac{1-2\varepsilon}{\sigma}.\]
  It follows that $\gamma_{\sigma,d} \geq (1-2\varepsilon)/\sigma$. Since this is true for any $\varepsilon > 0$, we conclude that $\gamma_{\sigma,d} \geq 1/\sigma$.
\end{proof}
It is now easy to disprove that
$\gamma_{\sigma,d}=\gamma_{\sigma,2}^{d-1}$ for large $d$.
Indeed,
 since $\gamma_{\sigma,2} < 1$~\cite{cs75}, then $\lim_{d\to \infty}
\gamma_{\sigma,2}^{d-1} = 0$. On the other hand, the previous lemma
asserts that $\gamma_{\sigma,d} \geq 1/\sigma$ for all~$d$, hence for
$d$ large enough, $\gamma_{\sigma,2}^{d-1} < \gamma_{\sigma,d}$.

In particular, for the case $\sigma=2$, Lueker~\cite{lueker03} proved that
  $\gamma_{2,2}\leq U$ for $U=0.826280$.
Thus, for all $d\geq 5$, we have the strict inequality
\[
\gamma_{2,2}^{d-1} \leq (0.826280)^{d-1} < 1/2 \leq \gamma_{2,d}.
\]

\section{Disproving $\gamma_{2,3}=\gamma_{2,2}^{2}$}\label{sec:complex}
\subsection{Diagonal common subsequence}
As already mentioned, the best known provable lower bound for
  $\gamma_{2,2}$ found so far is due to Lueker~\cite{lueker03}.
The starting point of Lueker's lower bound technique is a
  result by Alexander~\cite{alexander94} who
  related the expected length of the LCS of two random
  strings of the same length $n$, to the expected
  length of the LCS of two random strings whose lengths sum up to $2n$.
Below, we establish an analog of Alexander's result but for the case
  of $d$ randomly chosen sequences.

Let $C[j..k]$ denote the substring $C[j]C[j+1]\dots C[k]$ formed by
  all the characters between the $j$-th and \mbox{$k$-th} positions of $C$.
Given strings $A_1,\ldots, A_d$ of length at least $n$, we say
  that $B$ is an $n$-diagonal common subsequence of $A_1,\ldots, A_d$ if
  $B$ is a common subsequence of a set of prefixes of $A_1,\ldots, A_d$
  whose lengths sum to~$n$, i.e., if for some
  indices $i_1,\ldots, i_d$ such that $i_1 + \dots + i_d = n$, the
  string $B$ is a common subsequence of
  $A_1[1..i_1],\ldots,A_d[1..i_d]$.

Let $D_n(A_1,\ldots, A_d)$ denote the length of a longest
  $n$-diagonal common subsequence of the strings $A_1,\ldots, A_d$.
We denote by $\DD_{n,\sigma,d}$ the random variable $D_n(A_1,\ldots, A_d)$
  where the strings $A_1,\ldots, A_d$ are chosen according to
  $\calU_{n,\sigma}$.

The main objective of this section is to prove the following
  extension of a result of Alexander~\cite[Proposition~2.4]{alexander94}
  for the $d=2$ case:
\begin{theorem}\label{th:Alex2}
For all $n\geq d$,
\[ d\cdot \E\DD_{n,\sigma,d} - d^{3/2}\sqrt{2n\ln n} \leq \E\LL_{n,\sigma,d} \leq \E\DD_{nd,\sigma,d}.\]
In particular, for all $\sigma$ there exists $\delta_{\sigma,d}$ such that:
\[
\delta_{\sigma,d} = \lim_{n \to \infty} \frac{\E \DD_{n,\sigma,d}}{n} = \frac{\gamma_{\sigma,d}}{d}.
\]
\end{theorem}
For the sake of clarity of exposition, before proving Theorem~\ref{th:Alex2}
  we establish some intermediate results.

\begin{lemma}\label{lem:lower}
For all $n$ and $d$,
  $\E\LL_{n,\sigma,d} \leq \E\DD_{nd,\sigma,d}$.
\end{lemma}
\begin{proof}
Let $A_1,\ldots,A_d$ be random strings independently chosen
  according to $\calU_{nd,\sigma}$.
Since a longest common subsequence of $A_1[1..n],\ldots,A_d[1..n]$
  is also an $nd$-diagonal common subsequence of $A_1,\ldots, A_d$,
\[
L(A_1[1..n],\ldots,A_d[1..n]) \leq D_{nd}(A_1,\ldots,A_d).
\]
Taking expectation on both sides of the previous inequality yields
  the desired conclusion.
\end{proof}

\begin{lemma}\label{lem:upper}
For all $n\geq d$,
\[
d\cdot \E\DD_{n,\sigma,d} - d^{3/2}\sqrt{2n\ln n}
  \leq \E\LL_{n,\sigma,d}.
\]
\end{lemma}
\begin{proof}
Let $A_1,\ldots, A_d$ be a list of words of length $n$. Note that if
we change one character of any word in the list, then the values
$L(A_1,\ldots, A_d)$ and $D_n(A_1,\ldots, A_d)$ will change by at most
one unit. It follows that the random variables $\LL_{n,\sigma,d}$ and
$\DD_{n,\sigma,d}$ (seen as functions from $(\Sigma^{n})^{d}$ to $\RR$)
are both $1$-Lipschitz. Applying Azuma's inequality
  (as treated in for example~\cite[\S~2.4]{jlr00})
  we get:
\[
\Pr\left[\DD_{n,\sigma,d} \leq \E\DD_{n,\sigma,d} - \sqrt{n/2}\right]
  \leq \exp\left(-\frac{2(n/2)}{nd}\right) = e^{-1/d}
  < \frac{d}{d+1},
\]
where the last inequality holds since $e^{-x}< 1/(x+1)$ for
  all $x>0$.

Let $\lambda = \E\DD_{n,\sigma,d} - \sqrt{n/2}$.
Since $\DD_{n,\sigma,d} > \lambda$ implies that
  there are positive indices $i_1,\ldots,i_d$ such that
  $i_1 + \dots + i_d = n$ and
  $L(A_1[1..i_1],\ldots,A_d[1..i_d]) \geq \lambda$,
\[
  \Pr[\DD_{n,\sigma,d} > \lambda]
  \leq \sum_{\substack{0<i_1,\ldots,i_d<n,\\ i_1 + \ldots + i_d = n}}
  \Pr[L(A_1[1..i_1],\ldots,A_d[1..i_d]) > \lambda].
\]
Let $I$ be the number of summands in the right hand side.
Note that $I=\binom{n-1}{d-1}$ since it counts the
   the number of ways of partitioning $n$ into $d$ positive
  summands.
It follows that there exist positive
  $j_1,\ldots, j_d$ summing to $n$ such that:
\[
\Pr[L(A_1[1..j_1],\ldots,A_d[1..j_d]) > \lambda]
  >   \frac{1}{I}\left(1-\frac{d}{d+1}\right)
  = \frac{1}{I(d+1)}.
\]

Note that the distribution of
  the random variable $L(A_1[1..j_1], \ldots,A_d[1..j_d])$ is the same as
  the distribution of
  $L(A_1[1..j_{\tau(1)}],\ldots,A_d[1..j_{\tau(d)}])$ for any permutation $\tau: [d] \to [d]$.
It is also easy to see that the distribution of
  $L(A_1[a_1..b_1],\ldots,A_d[a_d..b_d])$ and
  $L(A_1[a'_1..b'_1],\ldots,A_d[a'_d..b'_d])$ are the same
  when $b_m-a_m=b'_m-a'_m$ for all $1\leq m \leq d$.

Now, let $\tau$ be the cyclic permutation
  $(1 2 \ldots d)$ and for $0\leq m\leq d-1$
  let $\calE_m$ denote the event
\[
  L\left(A_1\left[\sum_{l=0}^{m-1} j_{\tau^l(1)} + 1\ .. \sum_{l=0}^{m} j_{\tau^l(1)}\right],\ldots,A_d\left[\sum_{l=0}^{m-1} j_{\tau^l(d)} + 1\ .. \sum_{l=0}^{m} j_{\tau^l(d)}\right]\right) > \lambda.
\]
In particular, $\calE_0$ is the event $\{L(A_1[1..j_1],\ldots,A_d[1..j_d]) > \lambda\}$ whose probability was bounded above. Note that the events $\calE_0,\ldots,\calE_{d-1}$ are equiprobable.
Since each of the $\calE_m$'s depends on a
  different set of characters, they are independent.
Moreover, if  $\calE_0,\ldots,\calE_{d-1}$ simultaneously occur, then
  by concatenating the common subsequences of each block of characters
  we get that $L(A_1,\ldots, A_d) > d\lambda $.
Hence,
\begin{equation}
\left(\frac{1}{I(d+1)}\right)^d <
  \prod_{m=0}^{d-1}\Pr[\calE_m] = \Pr[\calE_0,\calE_1,\ldots,\calE_{d-1}]
  \leq \Pr[\LL_{n,\sigma,d} > d\lambda]. \label{eq:dcs1}
\end{equation}
Applying Azuma's inequality again, we have:
\begin{align}
\Pr\left[\LL_{n,\sigma,d} \geq \E\LL_{n,\sigma,d} + \sqrt{\frac{nd^2\ln(I(d+1))}{2}}\right] 
&\leq \left(\frac{1}{I(d+1)}\right)^d. \label{eq:dcs2}
\end{align}
Combining \eqref{eq:dcs1} and \eqref{eq:dcs2}
  and recalling that $\lambda = \E\DD_{n,\sigma,d} - \sqrt{n/2}$
  we obtain:
\[
\Pr\left[\LL_{n,\sigma,d}
  \geq \E\LL_{n,\sigma,d} + \sqrt{\frac{nd^2\ln(I(d+1))}{2}}\right]
  < \Pr\left[\LL_{n,\sigma,d} > d\E\DD_{n,\sigma,d} - d\sqrt{\frac{n}{2}}\right].
\]
Hence:
\[
\E\LL_{n,\sigma,d} + \sqrt{\frac{nd^2\ln(I(d+1))}{2}}
  \geq d\E\DD_{n,\sigma,d} - d\sqrt{\frac{n}{2}}.
\]
Since $2\leq d\leq n$,
  $(d+1)I = (d+1)\binom{n-1}{d-1} \leq n^d$, and so:
\[
  d\E\DD_{n,\sigma,d} \leq \E\LL_{n,\sigma,d}  + d\sqrt{\frac{n}{2}}
  + \sqrt{\frac{nd^2\ln(I(d+1))}{2}}
  \leq \E\LL_{n,\sigma,d} + d^{3/2}\sqrt{2n\ln(n)}.
\]
\end{proof}

\begin{prooff}{Theorem~\ref{th:Alex2}}
Lemmas~\ref{lem:lower} and~\ref{lem:upper} already give the bounds
  on $\E L_{n,\sigma,d}$.

To complete the proof we need to show that
  $\lim_{n \to \infty} \E \DD_{n,\sigma,d}/n$ exists and that
  its value is $\gamma_{\sigma,d}/d$.
By Lemmas~\ref{lem:lower} and~\ref{lem:upper} we have:
\[
\E\LL_{n,\sigma,d} \leq \E\DD_{nd,\sigma,d}
  \leq \frac{1}{d}\E\LL_{nd,\sigma,d} + d^{1/2}\sqrt{2nd \ln(nd)}.
\]
Dividing by $n$,
  it follows that $\lim_{n\to \infty}\E\DD_{nd,\sigma,d}/n=\gamma_{\sigma,d}$.
Furthermore, $\E\DD_{n,\sigma,d}$ is non decreasing in $n$, so:
\[
\frac{\lfloor n/d\rfloor}{n/d} \cdot \frac{\E\DD_{d\lfloor{n/d}\rfloor,\sigma,d}}{\lfloor{n/d}\rfloor}
  \leq \frac{\E\DD_{n,\sigma,d}}{n/d}
  \leq \frac{\lceil n/d\rceil}{n/d} \cdot \frac{\E\DD_{d\lceil{n/d}\rceil,\sigma,d}}{\lceil{n/d}\rceil}.
\]
Since both the left hand side and right hand side terms above
  converge to $\gamma_{\sigma,d}$ when $n\to\infty$,
  the middle term also converges to
  that value,
  and so $\lim_{n \to \infty} \E \DD_{n,\sigma,d}/n  = \gamma_{\sigma,d}/d$
  as claimed.
\end{prooff}

\subsection{Longest common subsequence of two words over a binary alphabet}
In this section we describe Lueker's~\cite{lueker03} approach
  for finding a lower bound on $\gamma_{d,\sigma}$ when $d=\sigma=2$.
Later on, we will generalize Lueker's technique to the cases of
  arbitrary $d$ and $\sigma$.

Let $X_1$ and $X_2$ be two random sequences chosen from $\calU_{n,2}$,
  i.e.~strings of length $n$ such that all their characters are
  chosen uniformly and independently from the binary alphabet $\{0,1\}$.
Lueker defines, for any two strings $A$ and $B$ over the binary
  alphabet, the quantity
\[
W_n(A,B) = \E\left[\max_{i+j=n}L(AX_1[1..i],BX_2[1..j])\right].
\]
Informally, $W_n(A,B)$ represents the expected length of a LCS of two
  strings with prefixes $A$ and $B$ respectively and suffixes
  formed by uniformly and independently choosing $n$ characters in
  $\{0,1\}$.
It is easy to see that
  $W_n(A,B)$ behaves as $\DD_{n,2,2}$ as
  $n\to\infty$.
Hence, applying Alexander's $d=2$ version of Theorem~\ref{th:Alex2},
  Lueker observes that for all $A,B\in\{0,1\}^{*}$,
\[
\gamma_{2,2} = \lim_{n\to\infty} \frac{W_{2n}(A,B)}{n}.
\]
A natural idea is to approximate $\gamma_{2,2}$ by $W_{2n}(A,B)/n$.
Fix the length $l \in \NN$ of the strings $A$ and $B$ and
  denote by $w_n$ the $2^{2l}$ dimensional vector
  whose coordinates correspond to the values $W_n(A,B)$ when
  $A$ and $B$ vary over all binary sequences of length $l$.
For example, when $l=2$ the vector $w_n$ has the following form:
\[
w_n = \begin{pmatrix}w_n[00,00]\\w_n[00,01]\\\vdots\\w_n[11,10]\\w_n[11,11]\end{pmatrix} = \begin{pmatrix}W_n(00,00)\\W_n(00,01)\\\vdots\\W_n(11,10)\\W_n(11,11)\end{pmatrix}.
\]
Lueker established a lower bound for each component of
$w_n$ as a function of the components of~$w_{n-1}$ and $w_{n-2}$. To
reproduce that lower bound, we need to introduce some more notation.
If $A = A[1]A[2]\ldots A[l]$ is a sequence of length $l\geq 2$, let
  $h(A)$ denote the \emph{head} of $A$, i.e.~its first character,
  and let $T(A)$ denote its \emph{tail}, i.e.~the substring obtained
  from $A$ by removing its head.
In other words, $h(A) = A[1]$ and $T(A)=A[2..l]$.
It is easy to see that the following relations among
  $w_n$,$w_{n-1}$ and $w_{n-2}$ hold:

\begin{itemize}
\item If $h(A)=h(B)$, then
\[ w_n[A,B] \geq
      1 + \frac{1}{4}\sum_{(c,c')\in \{0,1\}^2} w_{n-2}[T(A)c,T(B)c'].
\]

\item  If $h(A)\neq h(B)$, then
\[
    w_n[A,B] \geq
      \frac{1}{2} \max\left\{\sum_{c \in \{0,1\}}\! w_{n-1}[T(A)c,B],\sum_{c \in \{0,1\}}\! w_{n-1}[A,T(B)c]\right\}.
\]
\end{itemize}
Using the previous inequalities one can define a function $F: \RR^{2^{2l}} \times \RR^{2^{2l}} \to \RR^{2^{2l}}$ such that for all $n\geq 2$, we have $w_n\geq F(w_{n-1},w_{n-2})$. Furthermore, the function $F$ can be decomposed in two simpler functions $F_=$ and $F_{\neq}$ such that if $\Pi_=$ and $\Pi_{\neq}$ are the projections of the vectors onto the coordinates corresponding to the pairs of words with the same and different heads respectively, then:
\[
\Pi_=(w_n) \geq F_=(w_{n-2}),
\qquad \text{and} \qquad
\Pi_{\neq}(w_n) \geq F_{\neq}(w_{n-1}).
\]
It might be useful to see some examples of these transformations.
For instance, to obtain a lower bound of $w_n[001, 011]$, one considers:
\begin{eqnarray*}
\lefteqn {w_n[001,011] \geq F_=(w_{n-2})[001,011] } \\
&&= 1 + \frac{1}{4}\left(w_{n-2}[010,110] + w_{n-2}[010,111] + w_{n-2}[011,110] + w_{n-2}[011,111]\right).
\end{eqnarray*}
And to bound $w_n[001,111]$,
\begin{eqnarray*}
\lefteqn{w_n[001,111] \geq F_{\neq}(w_{n-1})[001,111]} \\
&&= \frac{1}{2}\max\left\{w_{n-1}[010,111]+ w_{n-1}[011,111], w_{n-1}[001,110]+ w_{n-1}[001,111] \right\}.
\end{eqnarray*}

\subsection{Longest common subsequence of $d$ words over general
  alphabets}\label{ssn-lcs-d-words}
In this section we extend Lueker's lower bound arguments as described in
  the previous section to the general case of $d$ strings
  whose characters are uniformly and independently chosen
  over an alphabet of size $\sigma$.

Let $X_1,\ldots, X_d$ be a collection of $d$ independent random
  strings chosen according to $\calU_{n,\sigma}$ and
  let $A_1,\ldots,A_d$ be a collection of $d$ finite sequences
  over the same alphabet.
We now consider:
\[
W_n(A_1,\ldots,A_d) =
  \E\left[\max_{i_1+\ldots+i_d= n}L(A_1X_1[1..i_1],\ldots, A_dX_d[1..i_d])\right].
\]
This quantity represents the expected length of a LCS of $d$ words
  with prefixes $A_1,\ldots, A_d$ respectively and $d$
  suffixes whose lengths sum up to $n$ and whose
  characters are uniformly and independently chosen
  in $\Sigma=\{1,\ldots,\sigma\}$.
Since $W_n(A_1,\ldots,A_d)$ and $\DD_{n,\sigma,d}$
  behave similarly as $n\to\infty$,
  Theorem~\ref{th:Alex2} implies that for all
  $A_1,\ldots, A_d$,
\begin{equation}\label{eq:limitwn}
\gamma_{\sigma,d} =
  \lim_{n\to\infty} \frac{W_{nd}(A_1,\ldots, A_d)}{n}.
\end{equation}

Just as in the $d=2$ case, fix $l\in\NN$ and denote by
  $w_n$ the $\sigma^{ld}$ dimensional vector
  whose coordinates are all the values
  of $W_{nd}(A_1,\ldots, A_d)$ when $A_1,\ldots,A_d$ vary
  over all sequences in~$\Sigma^{l}$.
We again seek a lower bound for $w_n$ as a function of
  vectors $w_m$, with $m< n$.

It is easy to see that if all
  the strings $A_1,\ldots,A_d$ start with the same character, then:
\[
w_n[A_1,\ldots,A_d]
  \geq 1 + \frac{1}{\left|\Sigma^d\right|}\sum_{\vec{c} \in \Sigma^d}w_{n-d}[T(A_1)c(1), T(A_2)c(2),\ldots, T(A_d)c(d)].
\]
Informally, the previous inequality asserts that if all the words
start with the same character then the expected length of the LCS of
all of them, allowing $n$ random extra characters, is at least 1 (the
first character) plus the average of the expected length of the LCS of
the words obtained by eliminating the first character and
``borrowing'' $d$ of the $n$ random characters.

If not all the words start with the same character, we can still find
  a lower bound, but to write it down we need to introduce some
  additional notation.
For any two sets $X$ and $Y$ we follow the standard convention
  of denoting by $Y^X$ the set of all mappings from $X$ to $Y$.
Also, for a $d$-tuple of strings $A=(A_1,\ldots,A_d)$ and $z \in\Sigma$
  we denote by $N_z(A)$ the set of indices $j \in \{1,\ldots,d\}$ such that
  $A_j$'s head is not equal to $z$,
  i.e.~to the set of string indices \emph{not}
  starting with~$z$.
For a mapping  $c: N_{z}(A)\to\Sigma$
  we define $\tau_z(A,c)$ as the
  the $d$-tuple of strings obtained from $A$ by replacing
  each string $A_i$ that does not start with $z$
  by the sequence obtained by eliminating its first character and
  adding the character $c(i)$ at its tail.
Formally, $\tau_z(A,c) = (A'_1,\ldots,A'_d)$
  where
\[
A'_i = \begin{cases}
  A_i, &\text{if $h(A_i)=z$,}\\
  T(A_i)c(i), &\text{if $h(A_i)\neq z$.}
\end{cases}
\]
A crucial fact is that for a $d$-tuple of strings $A$, if its
  coordinates do not all start with the same character, then
\[
w_n[A] \geq \max_{z \in \Sigma} \frac{1}{\left|\Sigma^{N_z(A)}\right|}
  \sum_{c \in \Sigma^{N_z(A)}} w_{n-|N_z(A)|}[\tau_z(A,c)].
\]
Informally, each term over which the maximum is taken
  corresponds to the expected length of the LCS of the strings
  one would obtain by disregarding all first characters of sequences
  not starting with $z$, and concatenating to the tail of
  these strings an element randomly chosen over the alphabet $\Sigma$.

For the sake of illustration, consider the following
  example of the derived inequalities when $\sigma=2$ and $d=4$:
\begin{align*}
w_n[001,011,101,001] \geq \max\bigg\{
&\frac{1}{2} \sum_{c \in \{0,1\}^{\{3\}}} w_{n-1}[001,011,01c(3),001],\\
&\frac{1}{2^3}\sum_{c \in \{0,1\}^{\{1,2,4\}}}w_{n-3}[01c(1),11c(2),101,01c(4)]\bigg\}.
\end{align*}
In the previous example only the third string over which
  $w_n$ is evaluated does not start with~$0$.
Hence, the first term over which the maximum is taken
  is the average of the values of $w_{n-1}$ evaluated at the
  two possible $4$-tuples of strings obtained from $A$
  by removing the initial $1$ from the third string and adding a
  $0$ or $1$ final character.
On the other hand, $w_n$ is evaluated at three strings that do not
  start with a $1$.
Hence, the second term over which the maximum is taken is the average of
  the values of $w_{n-3}$ over all the $4$-tuples of strings obtained from
  $A$ by removing all the initial~$0$'s and adding a $0$ or $1$
  final character to those same strings.

Expressing all the derived inequalities in vector form we
  have that there is a function
  $F:(\RR^{\sigma^{ld}})^{d} \to \RR^{\sigma^{ld}}$ such that
\begin{equation}
w_n \geq F(w_{n-1},w_{n-2},\ldots,w_{n-d}).\label{eq:recurrence}
\end{equation}
For the ensuing discussion it will be convenient
  to rewrite $F$ in an alternative way.
For each $z \in \Sigma$ we define the linear transformation
  $F_z: (\RR^{\sigma^{ld}})^{d} \to \RR^{\sigma^{ld}}$ such that
\begin{equation}
F_z(v_1,\ldots,v_d)[A] =  \begin{cases}
\displaystyle \frac{1}{\left|\Sigma^{N_z(A)}\right|}\sum_{c \in \Sigma^{N_z(A)}} v_{|N_z(A)|}[\tau_z(A,c)], &\text{if $|N_z(A)|\neq 0$,}\\
0, &\text{if $|N_z(A)|= 0$.}\end{cases}
\label{eq:defi_F}
\end{equation}
Then, if we let $b\in \RR^{\sigma^{ld}}$ be the vector with value $1$ in the
  coordinates associated to $d$-tuples of strings of length $l$ starting
  all with the same character and $0$ in the rest of the coordinates,
  $F$ can be expressed as:
\begin{equation}\label{eq:defi_F2}
F(v_1,\ldots,v_d) = b + \max_{z \in \Sigma}F_z(v_1,\ldots,v_d).
\end{equation}

\subsection{Finding a lower bound for $\gamma_{\sigma,d}$}
In the preceding section we established that for any $d$-tuple of strings
  $A=(A_1,\ldots,A_d)$, each of length $l$, we have
  $\gamma_{\sigma,d}~=~\lim_{n\to\infty} w_{nd}[A]/n$.
To lower bound this latter quantity one is tempted to try the following
  approach:  (1) For a fixed word length $l$, compute explicitly $w_0,\ldots,w_{d-1}$,
  and, (2) Define a new sequence of vectors $(v_n)_{n\in \NN}$
  as $v_i = w_i$ for $0\leq i \leq d-1$, and then iteratively
  define $v_n = F(v_{n-1},v_{n-2},\ldots, v_{n-d})$, for all $n\geq d$.
Since $F$ is monotone and by~\eqref{eq:recurrence},
  we have that $v_n\leq w_n$ for every $n\in \NN$.
It is natural to fix an arbitrary
  $d$-tuple of strings $A=(A_1,\ldots,A_d)$ and
  estimate a lower bound for~$\gamma_{\sigma,d}$
  by $\lim_{n\to\infty}v_{nd}[A]/n$ for large enough~$n$.

Unfortunately, for the approach discussed in the previous paragraph
  to work one would need to determine
  for which values of $n$ the quantity $v_{nd}[A]/n$
  is effectively a lower bound for~$\gamma_{\sigma,d}$.
Indeed, $v_{nd}[A]/n$ does not even need to be increasing and
  $w_{nd}[A]/n$ equals $\gamma_{\sigma,d}$ only in the limit
  when $n\to\infty$.
We will pursue a different approach that relies
  on the next lemma which is a generalization of an
  observation by Lueker~\cite{lueker03} for the $d=\sigma=2$ case.
\begin{lemma}\label{lem:F}
Let $\calF:(\RR^{\sigma^{ld}})^d\to \RR^{\sigma^{ld}}$
  be a transformation that satisfies the following properties:
\begin{enumerate}
  \item \textbf{Monotonicity:} If the inequality
  $(v_1,v_2,\ldots,v_d)\leq (w_1,w_2,\ldots,w_d)$ holds
  component-wise, then the inequality $\calF(v_1,v_2,\ldots,v_d)\leq
  \calF(w_1,w_2,\ldots,w_d)$ also holds component-wise.

  \item \textbf{Translation invariance:} Let $\ones$ be the vector of
  ones in $\RR^{\sigma^{ld}}$ and $\vec{\ones} = (\ones,\ldots,\ones)$
  be the vector of ones in $(\RR^{\sigma^{ld}})^d$. Then, for any
  $r\in \RR$ and for all $(v_1,v_2,\ldots,v_d) \in
  (\RR^{\sigma^{ld}})^d$,
  \[
  \calF((v_1,v_2,\ldots,v_d) + r\vec{\ones}) = \calF(v_1,\ldots,v_d) + r\ones.
  \]
  \item \textbf{Feasibility:} There exists a \emph{feasible triplet
  for $\calF$}, i.e.~a $(u,r,\varepsilon)$ with $u \in
  \RR^{\sigma^{ld}}$, $r \in \RR$, and $0 \leq \varepsilon \leq r$
  such that:
  \[
  \calF(u + (d{-}1)r\ones, \ldots, u+2r\ones, u+r\ones, u)
    \geq u + (dr - \varepsilon)\ones.
  \]
\end{enumerate}
Then, for any sequence $(v_n)_{n\in \NN}$ of vectors in $\RR^{\sigma^{ld}}$
  such that $v_n \geq \calF(v_{n-1},\ldots,v_{n-d})$ for all $n\geq d$,
  there exists a vector $u_0$ in $\RR^{\sigma^{ld}}$ such that for all
  $n\geq 0$,
\begin{equation}
v_n \geq u_0 + n(r-\varepsilon)\ones.\label{eq:lemma-invariance}
\end{equation}
\end{lemma}
\begin{proof}
Let $\calF$ be a transformation satisfying the hypothesis of the lemma and
  $(u,r,\varepsilon)$ a feasible triplet for $\calF$.
Let $(v_n)_{n\in \NN}$ be a sequence of vectors as in the lemma's
  statement and let $\alpha \in \RR$ be large enough so that
  for all $j\leq d-1$,
\[
v_j + \alpha \ones \geq u + j(r-\varepsilon)\ones.
\]
For example, set $\alpha$ to be the largest
  component of the vector
  $\max_{0\leq j\leq d-1}(u+j(r{-}\varepsilon)\ones-v_j)$.

Note that $u_0 = u - \alpha \ones$ satisfies~\eqref{eq:lemma-invariance}
  for all $n\leq d-1$.
We will prove by induction that this holds for all $n\in \NN$.
Suppose that~\eqref{eq:lemma-invariance} holds up to $n-1$.
Using the inductive hypothesis we have:
\begin{eqnarray*}
\lefteqn{(v_{n-1},\ldots, v_{n-d})} \\
  && \geq (u_0 + (n{-}1)(r{-}\varepsilon)\ones, \ldots, u_0 + (n{-}j)(r{-}\varepsilon)\ones, \ldots, u_0 + (n{-}d)(r{-}\varepsilon)\ones)\\
&&	=(u + (d{-}1)r\ones, \ldots, u + (d{-}j)r\ones + (j{-}1)\varepsilon\ones, \ldots, u + (d{-}1)\varepsilon\ones)\ +
\\ &&\qquad
  ((n{-}d)(r{-}\varepsilon) - (d{-}1)\varepsilon -\alpha)\vec{\ones}\\
&&\geq (u + (d{-}1)r\ones, \ldots, u + (d{-}j)r\ones, \ldots, u)\ +
  ((n{-}d)(r{-}\varepsilon) - (d{-}1)\varepsilon- \alpha)\vec{\ones}.
\end{eqnarray*}
Evaluating $\calF$
  at the terms on both sides of the previous inequality we get,
  by monotonicity and translation  invariance, that
   \begin{align*}
     v_n &\geq \calF(v_{n-1},\ldots, v_{n-d})\\ &\geq
     \calF(u + (d{-}1)r\ones, \ldots, u + (d{-}j)r\ones, \ldots, u)
   + ((n{-}d)(r{-}\varepsilon) - (d{-}1)\varepsilon -\alpha)\ones.
   \end{align*}
Since $(u,r,\varepsilon)$ is a feasible triplet, it follows that:
   \begin{align*}
     v_n &\geq u + (dr-\varepsilon)\ones + ((n{-}d)(r{-}\varepsilon) - (d{-}1)\varepsilon -\alpha)\ones
\\ &
   = u -\alpha \ones + n(r-\varepsilon)\ones = u_0 + n(r-\varepsilon)\ones.
   \end{align*}
   This completes the proof.
   \end{proof}

{From} $F$'s definition it easily follows that $F$ is monotone and
  invariant under translations.
If we find a feasible triplet $(u,r,\varepsilon)$ for $F$ then,
  by Lemma~\ref{lem:F}, we can conclude that the sequence of
  vectors $(w_n)_{n\in \NN}$ satisfy
  $w_n \geq u_0 + n(r-\varepsilon)\ones$ for all $n$.
It follows from~\eqref{eq:limitwn} that:
\[
\gamma_{\sigma,d} \geq d(r-\varepsilon).
\]
The key point we are trying to make is that in order
  to establish a good lower bound for $\gamma_{\sigma,d}$ one
  only needs to exhibit a good feasible triplet, namely one such that
  $(r-\varepsilon)$ is as large as possible.

Empirically, one observes that for any set of initial vectors
  $v_0, \ldots, v_{d-1}$, if one makes
  $v_{n+d} = F(v_{n+d-1},\ldots,v_{n})$ for all $n\in\NN$, then the
  sequence $(v_n)_{n \in \NN}$ is such that $v_n/n$ seems to
  converge to a vector with all its components taking the same value.
In fact, one observes that for large values of $n$ the vectors
  $v_n$ and $v_{n+1}$ differ
  essentially by a constant (independent of $n$) times the all ones vector.
Roughly, there exists a real value $r$ such that
  $v_{n+1}-v_{n}$ is approximately $r\ones$ for all large enough $n$.
Since, by definition $v_{n+d}=F(v_{n+d-1},\ldots,v_{n+1},v_n)$, this
  implies that
\[
F(v_n + (d{-}1)r\ones, v_n + (d{-}2)r\ones, \ldots, v_n + r\ones, v_n) \sim v_n + dr\ones.
\]
It follows that one possible approach
  to find a feasible triplet is to consider an
  $n$ large enough so that the difference between $v_{n}$ and $v_{n-1}$ is
  essentially a constant times the all ones vector.
Then, set $u=v_n$, and define $r$ as the maximum value such that
  $v_{n} - v_{n-1} \geq r\ones$ and $\varepsilon$
  as the minimum possible value such that the
  triplet $(u,r,\varepsilon)$ is feasible for $F$.
The following result validates the approach just described.
\begin{lemma}\label{lem:cont}
Let $\calF:(\RR^{\sigma^{ld}})^d\to \RR^{\sigma^{ld}}$ be a
  monotone and translation invariant transformation.
Let $v_0,\ldots,v_{d-1}\in\RR^{\sigma^{ld}}$ and
  $v_{n+d}=\calF(v_{n+d-1},\ldots,v_{n+1},v_n)$ for all $n\in\NN$.
If for some $r\in\RR$, $n_0\geq 1$ and $\varepsilon>0$ we have
  $||v_{n+1}-v_{n}-r\ones||_{\infty}\leq \varepsilon/2d$
  for all $n\in\{n_{0},\ldots,n_{0}{+}d{-}1\}$, then
  $(v_{n_0},r,\varepsilon)$ is a feasible triplet for $\calF$.
\end{lemma}
\begin{proof}
First, observe that the monotonicity and translation invariance
  property of $\calF$ implies that
\[
||\calF(x_0,\ldots,x_{d-1})-\calF(y_0,\ldots,y_{d-1})||_{\infty}
  \leq \max_{i=0,\ldots,d-1} ||x_{i}-y_{i}||_{\infty}\,.
\]
Let $u=v_{n_0}$ and note that $||v_{n_{0}+i}-(u+ir\ones)||_{\infty}
  \leq i\varepsilon/2d<\varepsilon/2$ for $0\leq i\leq d$.
Hence, by definition of $v_{n_0+d}$,
\[
||v_{n_0+d} -
  \calF(u+(d{-}1)r\ones,u+(d{-}2)r\ones, \ldots, u+r\ones,u)||_{\infty}
  \leq \varepsilon/2.
\]
Since $||v_{n_{0}+d}-(u+dr\ones)||_{\infty} \leq \varepsilon/2$ it
  follows that
\[
||(u+dr\ones)-
  \calF(u+(d{-}1)r\ones,u+(d{-}2)r\ones, \ldots, u+r\ones,u)||_{\infty}
  \leq \varepsilon.
\]
In other words, $(u,r,\varepsilon)$ is a feasible triplet for $\calF$.
\end{proof}
It is easy to check that $F$ satisfies the hypothesis of
  Lemma~\ref{lem:cont}.
This justifies, together with the empirical observation that
  $v_{n+1}-v_{n}$ is approximately $r\ones$ for large values
  of $n$, the general approach described in this section for
  finding a feasible triplet for $F$, and thus a lower bound
  for~$\gamma_{\sigma,d}$.
It is important to stress here that there is no
  need to prove the convergence of $v_n/n$ to~$r\ones$
  in order to establish the lower bound $\gamma_{\sigma,d}\geq d(r-\varepsilon)$.
We only need to find a feasible triplet $(u,r,\varepsilon)$ for $F$.
The characteristics of $F$, empirical observations and
  Lemma~\ref{lem:cont}, efficiently lead to such feasible
  triplets.

\subsection{Implementation and results. New bounds}
In this section we describe the procedure we implemented
  in order to find a feasible triplet $(u,r,\varepsilon)$ for $F$
  and, as a corollary, a lower bound for~$\gamma_{\sigma,d}$.
The procedure is called \textsc{FeasibleTriplet}, it is parameterized
  in terms of the number of sequences $d$ and the alphabet $\Sigma$,
  and its pseudocode is given in Algorithm~\ref{algo}.
\begin{algorithm}[ht]
\caption{Procedure for computing a feasible triple for $F$}
\label{algo}
\begin{algorithmic}[1]
\Procedure{FeasibleTriplet$_{d,\Sigma}$}{$l,n$}
\Comment{$l\in\NN$ parameter, $n\in\NN$ iteration steps}
\For{$i=0,\ldots,d-1$}
  \State $v_i\gets\zeros$
  \Comment{Where $\zeros$ denotes the vector of zeros in $\RR^{\sigma^{ld}}$}
\EndFor
\State $(u,r,\varepsilon) \gets (v_0,0,0)$
\For{$i=d,\ldots,n$}
  \State $v_i\gets F(v_{i-1}, v_{i-2},\ldots, v_{i-d})$
  \State $R \gets  \max_{A\in(\Sigma^{l})^{d}}{(v_{i}-v_{i-1})[A]}$
  \State $W\gets v_{i}+ dR\ones -
     F(v_{i}+(d{-}1)R\ones, \ldots, v_{i}+R\ones, v_{i})$
  \State $E\gets \max\{0, \max_{A\in(\Sigma^{l})^{d}}W[A]\}$
  \If{$R-E \geq r - \varepsilon$}
    \State $(u,r,\varepsilon) \gets (v_i,R,E)$
  \EndIf
\EndFor
\State \textbf{return} $(u,r,\varepsilon)$
\EndProcedure
\end{algorithmic}
\end{algorithm}
In order to implement $F$ we rely on the characterization
  given by \eqref{eq:defi_F} and \eqref{eq:defi_F2}.
Since the $F_z$'s are linear transformations, they can
  be represented as matrices.
This allows for fast evaluation of the $F_z$'s,
  but requires a prohibitively large amount of main memory
  for all but small values of $\sigma$, $l$ and $d$.
In order to optimize memory usage, we use the fact that
  by distinguishing~(\ref{eq:defi_F}) according to the cardinality of
  $N_{z}(A)$ where $A\in(\Sigma^{l})^{d}$,
  $F_z$ can be written as:
\[
F_z(v_1,\ldots, v_d) = \frac{1}{\sigma^1}F_{z,1}(v_1)
  + \ldots +\frac{1}{\sigma^{d}}F_{z,d}(v_d),
\]
where
\[
F_{z,i}(v_i)[A] =  \begin{cases}
\displaystyle \sum_{c \in \Sigma^{N_z(A)}}
  v_{i}[\tau_z(A,c)], &\text{if $|N_z(A)|=i$,}\\
0, &\text{otherwise.}
\end{cases}
\]
Note in particular that every $F_{z,i}$ can be represented as a
  $0$-$1$ sparse matrix.

In our experiments we ran Algorithm~\ref{algo} for different
  values of $l$ and alphabet sizes $\sigma$.
As one would expect, the derived lower bounds improve as $l$
  grows.
However, the memory resources required to perform the computation
  also increases.
Indeed, throughout the second loop of Algorithm~\ref{algo}
  we need to store $d$ vectors of dimension $\sigma^{ld}$.
Also, a simple analysis of the definition of the sparse matrix
  $F_{z,i}$ shows that it has
  $\binom{d}{i}\sigma^{(l-1)d}(\sigma-1)^{i}\sigma^{i}$ non-zero
  entries.
It follows that a sparse matrix representation of $F_{z}$ has
  roughly $\sigma^{ld}(\sigma-1)^{d}$ non-zero entries.
Hence, the necessary computations are feasible only for
  small values of $\sigma$, $l$ and $d$, unless additional
  features of the matrices involved are taken advantage of in
  order to optimize memory usage.

Table~\ref{table:results2} summarizes the lower bounds we obtain
  for $\gamma_{\sigma,2}$ and contrasts them with previously derived
  ones.
To the best of our knowledge, for the $d=2$ case and alphabet sizes
  $3$, $4$, $5$, and $6$, this work provides the currently best
  known lower bounds for $\gamma_{\sigma,2}$.
It might be worth mentioning that, as can be seen in that table, the
bound of~\cite{dancikthesis,deken79} is better than the bound
of the more recent work of~\cite{bgns99} for alphabet size $6$, and that for bigger alphabet
sizes, the bound of~\cite{dancikthesis,deken79} is still better than ours.

The best known lower bound for $\gamma_{2,2}$ is still the one
  established by Lueker~\cite{lueker03}.
Table~\ref{table:results} lists the distinct choices of $\sigma$
  and $d$ for which we could execute Algorithm~\ref{algo} and
  indicates the value of the parameter $l$ giving rise to the
  reported lower bound.

\begin{table}[!ht]
\centering
 \begin{tabular}{|c|c|c|c|}
 \hline 
           & \multicolumn{3}{|c|}{$\gamma_{\sigma,2}$}\\
  \cline{2-4}
  $\sigma$ & This work & Baeza et.~al. lower bound~\cite{bgns99} & Dan\v{c}\'{\i}k-Deken's lower bound~\cite{dancikthesis,deken79}\\
  \hline
  3        & {\bf 0.671697}        & 0.63376                  &0.61538\\
  4        & {\bf 0.599248}        & 0.55282                  &0.54545 \\
  5        & {\bf 0.539129}        & 0.50952                  &0.50615\\
  6        & {\bf 0.479452}        & 0.46695                  &0.47169\\
  7        & 0.444577              & -                        &{\bf 0.44502}\\
  8        & 0.356545              & -                        &{\bf 0.42237}\\
  9        & 0.327935              & -                        &{\bf 0.40321}\\
  10       & 0.303490              & -                        &{\bf 0.38656}\\
  \hline
 \end{tabular}
\caption{Best known lower bounds for $\gamma_{\sigma,2}$ (in boldface). }\label{table:results2}
\end{table}

\begin{table}[!ht]
\begin{tabular}{cc}
    \begin{minipage}{0.45\linewidth}
    \centering
        \begin{tabular}{|ccc|}
        \hline 
        \multicolumn{ 3}{|c|}{Alphabet size $\sigma=2$} \\
        \hline
        $d$ & $L$ such that $\gamma_{2,d}\geq L$ & Parameter $l$  \\
        \hline
         2 &   0.781281 &          10\\
         3 &   0.704473 &          7 \\
         4 &   0.661274 &          5 \\
         5 &   0.636022 &          4 \\
         6 &   0.617761 &          3 \\
         7 &   0.602493 &          2 \\
         8 &   0.594016 &          2 \\
         9 &   0.587900 &          2 \\
        10 &   0.570155 &          1 \\
        11 &   0.570155 &          1 \\
        12 &   0.563566 &          1 \\
        13 &   0.563566 &          1 \\
        14 &   0.558494 &          1 \\
        \hline 
        \end{tabular}
    \vspace{8pt}

        \begin{tabular}{|ccc|}
        \hline 
        \multicolumn{ 3}{|c|}{Alphabet size $\sigma=3$} \\
        \hline
        $d$ & $L$ such that $\gamma_{3,d}\geq L$ & Parameter $l$  \\
        \hline
         2 &   0.671697 &          6 \\
         3 &   0.556649 &          4 \\
         4 &   0.498525 &          3 \\
         5 &   0.461402 &          2 \\
         6 &   0.421436 &          1 \\
         7 &   0.413611 &          1 \\
         8 &   0.405539 &          1 \\
        \hline
        \end{tabular}
    \vspace{8pt}

        \begin{tabular}{|ccc|}
        \hline
        \multicolumn{ 3}{|c|}{Alphabet size $\sigma=4$} \\
        \hline
        $d$ & $L$ such that $\gamma_{4,d}\geq L$ & Parameter $l$  \\
        \hline
         2 &   0.599248 &          5 \\
         3 &   0.457311 &          3 \\
         4 &   0.389008 &          2 \\
         5 &   0.335517 &          1 \\
         6 &   0.324014 &          1 \\
        \hline
        \end{tabular}
    \end{minipage}
&
    \begin{minipage}{0.45\linewidth}
    \centering
        \begin{tabular}{|ccc|}
        \hline
        \multicolumn{ 3}{|c|}{Alphabet size $\sigma=5$} \\
        \hline
        $d$ & $L$ such that $\gamma_{5,d}\geq L$ & Parameter $l$  \\
        \hline
         2 &   0.539129 &          4 \\
         3 &   0.356717 &          2 \\
         4 &   0.289398 &          1 \\
         5 &   0.273884 &          1 \\
        \hline
        \end{tabular}
    \vspace{8pt}

        \begin{tabular}{|ccc|}
        \hline
        \multicolumn{ 3}{|c|}{Alphabet size $\sigma=6$} \\
        \hline
        $d$ & $L$ such that $\gamma_{6,d}\geq L$ & Parameter $l$  \\
        \hline
         2 & 0.479452 &          3 \\
         3 & 0.309424 &          2 \\
         4 & 0.245283 &          1 \\
        \hline
        \end{tabular}
    \vspace{8pt}

        \begin{tabular}{|ccc|}
        \hline
        \multicolumn{ 3}{|c|}{Alphabet size $\sigma=7$} \\
        \hline
        $d$ & $L$ such that $\gamma_{7,d}\geq L$ & Parameter $l$  \\
        \hline
         2 & 0.444577 &          3 \\
         3 & 0.234567 &          1 \\
         4 & 0.212786 &          1 \\
         \hline
        \end{tabular}
    \vspace{8pt}

        \begin{tabular}{|ccc|}
        \hline
        \multicolumn{ 3}{|c|}{Alphabet size $\sigma=8$} \\
        \hline
        $d$ & $L$ such that $\gamma_{8,d}\geq L$ & Parameter $l$  \\
        \hline
         2 & 0.356545 &          2 \\
         3 & 0.207547 &          1 \\
        \hline
        \end{tabular}
    \vspace{8pt}

        \begin{tabular}{|ccc|}
        \hline
        \multicolumn{ 3}{|c|}{Alphabet size $\sigma=9$} \\
        \hline
        $d$ & $L$ such that $\gamma_{9,d}\geq L$ & Parameter $l$  \\
        \hline
         2 &        0.327935 &          2 \\
         3 &        0.186104 &          1 \\
        \hline
        \end{tabular}
    \vspace{8pt}

        \begin{tabular}{|ccc|}
        \hline
        \multicolumn{ 3}{|c|}{Alphabet size $\sigma=10$} \\
        \hline
        $d$ & $L$ such that $\gamma_{10,d}\geq L$ & Parameter $l$  \\
        \hline
         2 & 0.303490 &          2 \\
         3 & 0.168674 &          1 \\
        \hline
        \end{tabular}
    \end{minipage}
\end{tabular}
\caption{Lower bounds for $\gamma_{\sigma,d}$}\label{table:results}
\end{table}

\subsection{Disproving Steele's $\gamma_{2,2}=\gamma_{2,3}^2$ speculation}
We showed in Section~\ref{sec:simple} that
  $\gamma_{2,d} > \gamma_{2,2}^{d-1}$ for all $d\geq 5$.
We now establish that this is also the case when $d=3$ and $d=4$.
Recall that Lueker~\cite{lueker03} proved that
  $\gamma_{2,2}\leq U$ for $U=0.826280$.
{From} Table~\ref{table:results} we see that for $d=3$ and $d=4$,
  the indicated lower bound for $\gamma_{2,d}$ is strictly
  greater than~$U^{d-1}$, and therefore, is also strictly greater than
  $\gamma_{2,2}^{d-1}$.
This implies that   $\gamma_{2,d} > \gamma_{2,2}^{d-1}$ for
  $d=4$ and $d=3$ as claimed. Together with the results of Section~\ref{sec:simple} this establishes that $\gamma_{2,d} > \gamma_{2,2}^{d-1}$ for all $d\geq 3$.

\section{Final comments}\label{sec:conclusion}
As already mentioned at the start of this paper, Steele~\cite{steele86}
  pointed out that it would be of interest to find relations
  between the values of the $\gamma_{\sigma,d}$'s, especially
  between $\gamma_{2,2}$ and $\gamma_{2,3}$.
We think it would be very interesting if such a relation would
  exist.
In fact, it might shed some light upon the longstanding open
  problem of determining the exact value of the Chv\'atal-Sankoff
  constant.

Lacking a relation among the $\gamma_{\sigma,d}$'s it
  would still be interesting to relate these terms to
  some other constants that arise in connection with other
  combinatorial problems.
A step in this direction was taken by Kiwi, Loebl and
  Matou\v{s}ek~\cite{klm05} who showed that
  $\sqrt{\sigma}\gamma_{\sigma,2}\to c_2$
  when $\sigma\to\infty$, where $c_{2}$ is a constant that
  turns up in the study of the Longest Increasing Sequence (LIS)
  problem (also known as Ulam's problem).
Specifically, $c_2$ is the limit to which the expected
  length of a LIS of a randomly chosen permutation of $\{1,\ldots,n\}$
  converges when normalized by $\sqrt{n}$.
Logan and Shepp~\cite{ls77} and
  Vershik and Kerov~\cite{vk77} showed that $c_{2}=2$.
Consider now the following experiment: Choose $n$ points
  in a unit $d$-dimensional cube $[0,1]^{d}$ and let $H_{d}(n)$ be
  the random variable corresponding to the
  length of a longest chain (for the standard partial order
  in $\RR^{d}$) of the $n$ chosen points.
Bollob\'as and Winkler~\cite{bw88} proved that there are constants
  $c'_2, c'_3,\ldots$ such that $c'_d<e$,
  $\lim_{d\to\infty} c'_d=e$ and
  $\lim_{n\to\infty} H_{d}(n)/n^{1/d}=c'_d$.
By labeling a set $S$ of points in $[0,1]^{2}$ in
  increasing order of their $x$ coordinate and reading the
  labels in the order of their $y$ coordinates one can
  associate a permutation $\pi$ to the set $S$.
It is easy to see that a chain of points in $S$
  is in one to one correspondence to an increasing
  sequence of $\pi$.
Hence, it follows that $c'_2=c_2$.
Soto~\cite{sotothesis} extended the results of~\cite{klm05} and
  showed that $\sigma^{1-1/d}\gamma_{\sigma,d}\to c'_d$ when
  $\sigma\to\infty$.
We think that any similar type of result, or even a reasonable
  conjecture, that would hold for fixed $\sigma$ and $d$ would
  also be quite interesting.

\bibliographystyle{alpha}
\bibliography{biblio}

\end{document}